\documentclass[11pt]{article}
\usepackage{amssymb,amsmath,amsthm,fullpage,mathtools,enumitem, mathrsfs,mathtools}
\usepackage{tikz}
\usetikzlibrary{positioning,chains,fit,shapes,calc}
\usepackage{xcolor}
\usepackage{faktor}
\usetikzlibrary{patterns,decorations}
\usepackage{ytableau}

\definecolor{greenblue}{rgb}{0,.5,.44}
\definecolor{redgreen}{rgb}{.67,.43,0}
\definecolor{redblue}{rgb}{.5,0,.5}

\newtheorem{theorem}{Theorem}

\newtheorem{lemma}{Lemma}

\theoremstyle{definition}
\newtheorem{remark}{Remark}

\newtheorem*{definition}{Definition}
\newtheorem*{notation}{Notation}
\newtheorem*{definitions}{Definitions}
\newtheorem*{ack}{Acknowledgement}
\newtheorem{example}{Example}

\DeclarePairedDelimiterX\makeset[2]{\{}{\}}{#1\; \delimsize\vert\; #2}

\newcommand\off{
  \mathchoice
    {{\scriptstyle\mathcal{O}}}% \displaystyle
    {{\scriptstyle\mathcal{O}}}% \textstyle
    {{\scriptscriptstyle\mathcal{O}}}% \scriptstyle
    {\scalebox{.7}{$\scriptscriptstyle\mathcal{O}$}}%\scriptscriptstyle
  }

\newcommand{\aux}{\mathcal{A}}
\renewcommand{\P}{\mathcal{P}}
\newcommand{\C}{\mathcal{C}}
\newcommand{\M}{\mathcal{M}}
\newcommand{\Q}{\mathcal{Q}}
\newcommand{\I}{\mathcal{I}}

\DeclareMathOperator{\sP}{\mathscr{P}}
\DeclareMathOperator{\sB}{\mathscr{B}}

\DeclareMathOperator{\pr}{Pr}

\DeclareMathOperator{\seg}{seg}
\DeclareMathOperator{\cl}{Cl}

\newcommand{\stair}{\mathbb{S}}
\newcommand{\join}{\vee}
\newcommand{\myand}{\quad \text{and}\quad }

\begin{document}
\title{Rook and Wilf equivalence of integer partitions}

\author{
  Jonathan Bloom\\
  \texttt{Lafayette College}\\
  \texttt{bloomjs@lafayette.edu}\\
  \and 
  Dan Saracino\\
    \texttt{Colgate University}\\
    \texttt{dsaracino@colgate.edu}
}

\maketitle
\begin{abstract} The subjects of rook equivalence and Wilf equivalence have both attracted considerable attention over the last half-century.  In this paper we introduce a new notion of Wilf equivalence for integer partitions, and, using this notion, we prove that rook equivalence implies Wilf equivalence. We also prove that if we refine the notions of rook and Wilf equivalence in a natural way, then these two notions coincide.  In \cite{Bloom:On-cr2017} we prove that Wilf equivalence implies rook equivalence.
\end{abstract}

\section{Introduction}

In the late 1940's, Kaplansky and Riordon~\cite{kaplansky1946problem} introduced rook polynomials as a vehicle for the systematic study of permutations avoiding certain forbidden positions. (For example, derangements avoid the positions $(i,i)$ when considered as rook placements.)   As interesting objects of study in their own right, rook polynomials have since then received a great deal of attention (see,~\cite{barrese2014m,barrese2016bijections,briggs2006m,goldman1975rook}).   Of particular relevance to our work is the paper of Foata and Sch\"{u}tzenberger~\cite{foata1970rook} where the notion of rook equivalence of integer partitions was introduced and completely characterized.  

Around the same time, a parallel story unfolded in connection with the study of permutations avoiding certain forbidden patterns. Although originally touched on by MacMahon~\cite{MacMahon:Combinat15} in 1915,  this study was revived by Knuth in 1968 when he used permutations avoiding the pattern $231$ to characterize those permutations that are stack sortable~\cite{Knuth:The-art-68}.    This development launched a systematic study of permutation pattern avoidance. (Two books on the topic are~\cite{Bona:Combinat12} and \cite{kitaev2011patterns}.) The idea of Wilf equivalent patterns, i.e.,  patterns that are equally difficult to avoid, has played a fundamental role throughout this study.  

In recent years, the idea of ``pattern" avoidance, as it pertains to combinatorial objects other than permutations, has received growing attention.  Among some of the objects studied in this context have been words, set partitions, matchings, and Catalan objects (see~\cite{AlbertBouvel:A-genera2014,  BloomElizalde:Pattern-2013, BloomSaracino:Pattern-2016, Burstein:Enumerat1998,KlazarOn-abab-1996}).   In this paper we consider a new definition of pattern avoidance in the context of integer partitions and take up its systematic study.  In particular, we establish a connection between rook theory and pattern avoiding integer partitions by proving that rook equivalence implies Wilf equivalence of integer partitions.  We  also establish a partial converse by considering natural refinements of the notions of Wilf and rook equivalence.  It should be noted that Remmel's work unifying the theory of classical integer partition identities in~\cite{Remmel:Bijectiv1982} may also be viewed as a study of a  different definition of integer partition avoidance. We  will say more about this once exact definitions are established below. 

We begin by establishing some definitions and notation.  We first define $\P$ to be the set of all integer partitions and refine this for any $n\geq 0$ by letting $\P_n$ be the set of all integer partitions $\mu$ whose \emph{weight} (i.e., the sum of the parts) is $n$. We denote the weight of a partition $\mu$ by $|\mu|$.    It will be convenient in the sequel to adopt the convention that any integer partition $\mu$ consists of an infinite number of parts together with a part of size zero at infinity.  In other words, the parts of $\mu$ are such that $\mu_i\geq \mu_{i+1}\geq 0$ for each $i\in [1,\infty)$, only a finite number of the $\mu_i$  are positive, and $\mu_\infty =0$.    When writing partitions we only list the positive parts.  So the partition $(3,3,2,1,1,0,0,\ldots, 0)$ is, as is customary, simply written as $(3,3,2,1,1)$. This definition allows for the partition $(0,0,\ldots, 0)$ which we call the \emph{empty} partition.   

For any partition $\mu$ we define its \emph{height} to be the number of positive parts $h_\mu$.  Likewise, we define its \emph{width}, which we denote by $w_\mu$, to be the size of its first part $\mu_1$.  The empty partition is the only partition with height or width equal to 0.  Using these notions we refine our set $\P$ by defining, for any $h,k\geq 0$, the set
$$\P(h,k) = \makeset{\mu\in \P}{h_\mu \leq h, w_\mu = k}.$$

Although partitions are defined as weakly decreasing sequences of nonnegative integers, it is also common to view an integer partition as a Ferrers board.  For example, the partition $(3,3,2,1,1)$ can be drawn as the Ferrers board 
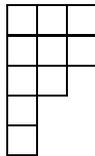
\ytableausetup{boxsize=1em}
$$\ydiagram{3,3,2,1,1}$$
consisting of 5 rows and 3 columns.  Viewed this way, the height and width of a partition become natural definitions.  Going forward we implicitly identify a partition with its Ferrers board.  As such any reference to a ``row" or ``column" of a partition should not cause confusion.  

Our next two definitions make immediate use of this identification.  First, for any partition $\sigma$ we denote by $\sigma^*$ the partition obtained by interchanging the rows and columns of (the Ferrers board for) $\sigma$.  We call $\sigma^*$ the \emph{conjugate} of $\sigma$.  For example, the conjugate of the above partition $(3,3,2,1,1)$ is 
$$(5,3,2) = \ydiagram{5,3,2}\ .$$
In what follows we make use of the fact that $\sigma_i^*$ is the length of the $i$th column of $\sigma$.    From this, it now follows easily that
$$\sigma_i^* = |\makeset{j}{\sigma_j\geq i}|.$$

The second definition that follows naturally from this partition/Ferrers board identification is new and at the heart of this paper.  We say a partition $\alpha$ \emph{contains} a partition $\mu$  provided that it is possible to delete rows and columns from $\alpha$ so that one obtains $\mu$.  For example, the partition 

$$(5,5,2,2,2) = \ydiagram{5,5,2,2,2}*[*(red)]{1+1,5,1+1,2,2}*[*(red)]{3+2}$$ 
contains $(2,1) = \ydiagram{2,1}$, since deleting the colored rows and columns yields $(2,1)$. In fact, it is easily seen that the only partition that does not contain $(2,1)$ is one whose Ferrers board is a rectangle. For another example, if $\mu=(4,3,3, 2,2)$ then by deleting the indicated rows and columns we see that $\mu$ is contained in the partition 
$$\ydiagram{6,5,5,5,4,4,2,2}*[*(red)]{0,5,0,0,4,0,2}*[*(red)]{2+2,2+2,2+2,2+2,2+2,2+2}.$$
  
For any partition $\mu$ we define $\P(\mu)$ to be the set of all partitions that contain $\mu$ and we define $\P_n(\mu) = \P(\mu) \cap \P_n$.  We say that $\mu,\tau\in \P$ are \emph{Wilf equivalent} provided that 
$$|\P_n(\mu)| = |\P_n(\tau)|$$
for all $n\geq 0$.  Similarly, we say that $\mu$ and $\tau$ are \emph{width-Wilf equivalent} provided that there are the same number of partitions of each weight and width that contain $\mu$ as there are that contain $\tau$.  

In what follows a further refinement of $\P(\mu)$ and $\P_n(\mu)$ is needed.    For any $k\geq 0$ we set
$$ \P(\mu,k) =\makeset{\alpha \in \P(\mu)}{w_\alpha=k+w_\mu}\ \textrm{and}\ \P_n(\mu,k) = \makeset{\alpha \in \P_n(\mu)}{w_\alpha = k+w_\mu}.$$
Note that $\mu$ and $\tau$ are width-Wilf equivalent if and only if $w_\mu = w_\tau$ and 
$$|\P_n(\mu,k)| = |\P_n(\tau,k)|$$
for all $n,k\geq 0$.

With our definition of integer partition containment established, we now compare it to Remmel's definition of partition containment as given in his landmark paper~\cite{Remmel:Bijectiv1982}.  In this work, Remmel stipulated that a partition $\alpha$ contains another partition $\mu$ if one can obtain $\mu$ by deleting only rows from $\alpha$.  Under this definition it is no longer the case that $(6,5,5,5,4,4,2,2)$ contains $(4,3,3,2,2)$, and it is easy to see that two partitions are Wilf equivalent if and only if they have the same weight. 

Remmel used his notion of containment in order to provide a unifying theory for the previously ad hoc study of partition identities.  For example, it has been well known since the time of Euler that the set $\mathcal{D}_n$ of partitions of weight $n$ that have distinct parts is equinumerous to the set $\mathcal{O}_n$  of partitions of weight $n$ that have only odd parts.   In terms of Remmel's language, the former set can be described as those partitions that do not contain any  member of 
$$\{(1,1), (2,2),(3,3),\ldots \},$$
while the latter set can be described as those partitions that do not contain any member of 
$$\{(2),(4),(6),\ldots\}.$$
In this way, the classic result that $|\mathcal{D}_n| = |\mathcal{O}_n|$ may be viewed as a result about a notion of Wilf equivalence.  Furthermore, one of the main results in~\cite{Remmel:Bijectiv1982} is a condition on subsets  $S,T\subseteq \P$  that guarantees Wilf equivalence.

The first goal of our paper is to provide a combinatorial description of the generating function
$$F_{\mu,k}(q) = \sum_{n\geq 0} |\P_n(\mu,k)|q^n$$ for $\P(\mu,k)$ 
directly in terms of the multiset of integers 
$$\{\mu_1 + 1, \mu_2+2, \mu_3+3, \ldots\}$$
and the parameter $k$ without appealing to the set of all partitions that contain $\mu$.  We prove this result in Section~\ref{sec:GF for P(mu,k)} modulo a few results whose proofs are deferred to Sections~\ref{sec:profiles} and~\ref{sec:staircases and augmented structures}. (These include Theorem~\ref{thm:Closures, Profiles, and Maxima}, which may be of interest in its own right.) In Section~\ref{sec:equivalence} we use this description of $F_{\mu,k}(q)$ to prove that rook equivalence implies Wilf equivalence (Theorem~\ref{thm:rook and Wilf}).  We also prove in this section (Theorem~\ref{thm:width Wilf}) that two partitions $\mu$ and $\tau$ with the same width are width-Wilf equivalent  if and only if they are rook equivalent.  Lastly, we prove (also in Theorem~\ref{thm:width Wilf}) the curious fact that if $w_\mu=w_\tau$  and
$$|\P_n(\mu,1)| = |\P_n(\tau,1)|$$
for all $n$, then $\mu$ and $\tau$ are actually width-Wilf equivalent. 

In \cite{Bloom:On-cr2017} we use Theorem~\ref{thm:rook and Wilf} to prove that Wilf equivalence implies rook equivalence, so these two notions coincide.

\section{The generating function $F_{\mu,k}$}\label{sec:GF for P(mu,k)}

The goal of this section is to prove Theorem~\ref{thm:GF of F_{mu,k}}, which gives a description of the generating function for $\P(\mu,k)$ directly in terms of the Ferrers board of $\mu$ and the parameter $k$ without appealing to all partitions that contain $\mu$.  The argument needed to arrive at this final description requires several steps.  The complete argument is presented here modulo a few requisite results whose proofs are postponed to Section~\ref{sec:profiles} and Section~\ref{sec:staircases and augmented structures}.  For clarity, we indicate such results when we state them here by enclosing them in boxes.   

We begin with some basic definitions.  
\begin{definition}
For any partitions $\mu$ and $\alpha$ we define their \emph{sum} to be the partition
$$\mu+\alpha= (\mu_1+\alpha_1, \mu_2+\alpha_2,\ldots).$$ 
\end{definition}

\ytableausetup{boxsize=.9em}	

For example, if $\mu = (2,2,2,1,1)$ and $\alpha = (3,2,2,1)$, then $\mu+\alpha= (5,4,4,2,1)$.  Graphically, the sum of $\mu$ and $\alpha$  is the partition whose columns are those of $\mu$ together with those of $\alpha$.  

\begin{definition}
For any partition $\beta$ we define $\Q(\beta)$ to be the set of all partitions obtained by  adding to $\beta$ an arbitrary number of parts no larger than $w_\beta$.  
\end{definition}
	 The generating function for $\Q(\beta)$ is easily seen to be
$$\sum_{\sigma\in \Q(\beta)} q^{|\sigma|} =\frac{q^{|\beta|}}{(1-q)(1-q^2)\cdots (1-q^{w_{\beta}})}.$$

Our first lemma allows us to rewrite the set $\P(\mu,k)$ as a union of these $\Q$-sets.

\begin{lemma}\label{lem:P(mu,k) is union} Fix $h,k > 0$ and a nonempty partition $\mu$ such that $h_\mu \leq h$.  Then
$$\P(\mu,k) = \bigcup_{\alpha\in \P( h,k)} \Q(\mu+\alpha).$$	
\end{lemma}

\begin{proof}
Clearly the righthand side is a subset of the lefthand side.  To see the reverse inclusion, fix a partition $\sigma\in \P(\mu,k)$ and identify a set $R$ of rows and a set $C$ of columns that may be deleted to obtain $\mu$.   

First, assume that the leftmost column  is not in $C$ and the top row is not in $R$.  As the leftmost column is not deleted, it follows that $R$ contains  exactly $h_{\sigma}-h_{\mu}$ rows.  Moreover, if we first delete these rows, then the fact that the top row remains means that we are left with a partition containing $\mu$ with height $h_\mu$ and width $w_\mu+k$.  Such a partition is clearly an element of the righthand side, since $h\geq h_{\mu}$.

It only remains to show that we can choose $R$ and $C$  that do not involve the top row or the leftmost column.    To see that this is possible, start with any choice of $R$ and $C$. If $1\in C$, then let $i$ be the leftmost column not deleted. (Such an $i$ exists since $\mu\neq \emptyset$.) Observe that we can now recover $\mu$ by deleting the columns in $C'=(C\cup \{i\})\setminus \{1\}$ together with all the rows in $R'=R\cup \{\sigma^*_i+1,\sigma^*_i+2,\ldots, h_{\sigma}\}$. (The fact that $\sigma^*_1>0$ guarantees that  $1\in R$ if and only if $1\in R'$.)
An analogous argument (applied to $\sigma^*$) shows that we can also modify our sets so that $1\notin R'$.  Note that the above parenthetical statement guarantees that applying the procedure twice does not unintentionally add the leftmost column back to $C'$.  This completes our proof.  
\end{proof}

It is important to note that the sets $\Q(\mu+\alpha)$ do not partition $\P(\mu,k)$.  For example, if $\mu$ and $\alpha$ are as in the above example and we let $\beta= (3,2,1,1)$, then
$$(5,4,4,3,2,1) \in \Q(\mu+\alpha)\cap \Q(\mu+\beta).$$  
This occurs because one could either start with $\mu+\alpha= (5,4,4,2,1)$ and add a part of size 3 or one could start with $\mu +\beta= (5,4,3,2,1)$ and add a part of size 4.  Consequently, we are forced to invoke the principle of Inclusion--Exclusion in order to use Lemma~\ref{lem:P(mu,k) is union} for enumerative purposes. For this to be an effective strategy, we need a concrete description of the sets
\begin{equation}\label{eq:intersection of Q's}
\bigcap_{\alpha \in P} \Q(\mu+\alpha),	
\end{equation}
where $P \subseteq \P(h, k)$. This leads us to our next definition.  

\begin{definition}
For any partitions $\alpha = (p^{m_p}, \ldots, 1^{m_1})$ and $\beta = (p^{n_p},\ldots, 1^{n_1})$, where $p$ is the largest integer that is a part of either $\alpha$ or $\beta$ and  $m_i,n_i\geq 0$, we define
$$\alpha \join \beta = \left(p^{e_p}, \ldots, 1^{e_1}\right),$$ 
where $e_i=\max(m_i,n_i)$.  For any nonempty finite set of partitions $P=\{\alpha^{(1)},\alpha^{(2)},\ldots, \alpha^{(s)}\}$ we define 
$$\vee_P(\mu):=\left(\mu + \alpha^{(1)}\right) \vee \cdots \vee \left(\mu + \alpha^{(s)}\right).$$	
\end{definition}

To see an example, take $P=\{\alpha,\beta\}$ where $\alpha$ and $\beta$ are as above.  Then
$$\vee_P(\mu)=(\mu+\alpha)\join (\mu+\beta) = (5,4,4,2,1) \vee (5,4,3,2,1) = (5,4,4,3,2,1).$$

The proof of our next lemma is left to the reader.    

\begin{lemma}\label{lem:intersection is multiset union}
Fix $h,k > 0$.  For any nonempty $P\subseteq \P(h,k)$ we have
	$$\bigcap_{\alpha\in P} \Q(\mu+\alpha) = \Q\left(\join_P(\mu)\right).$$

\end{lemma}
Using Lemma~\ref{lem:P(mu,k) is union} and Lemma~\ref{lem:intersection is multiset union}, we can take the first step toward determining $F_{\mu,k}(q)$. We have, for any nonempty $\mu$ such that $h_{\mu}\leq h$,
\begin{align}\label{eq:F1}
F_{\mu,k}(q) &=\sum_{\sigma\in \P(\mu,k)} q^{|\sigma|}\nonumber\\
 &=\sum_{\substack{P\subseteq \P(h,k)\\P\neq \emptyset}} (-1)^{|P|+1} \sum_{\sigma } q^{|\sigma|},\nonumber\\
 \shortintertext{where the inner sum is over all $\displaystyle\sigma\in\bigcap_{\alpha\in P}\Q(\mu+\alpha)$,   }
 &=\sum_{\substack{P\subseteq \P(h,k)\\P\neq \emptyset}} (-1)^{|P|+1} \sum_{\sigma\in \Q(\vee_P(\mu)) } q^{|\sigma|}\nonumber\\
 &=\sum_{\substack{P\subseteq \P(h,k)\\P\neq \emptyset}} (-1)^{|P|+1}  \frac{q^{|\join_P(\mu)|}}{(1-q)\cdots (1-q^{k+w_\mu})}.
\end{align}

Each of the partitions $\join_P(\mu)$ in this last sum is the value at $\mu$ of a function 
$\vee_P:\P\to\P$. To proceed with the determination of $F_{\mu,k}$ we will determine the conditions under which $\vee_{P} = \vee_{Q}$, for nonempty finite subsets $P,Q$ of $\P$.   

\begin{definitions}
For any partition $\sigma$ and $p\geq 0$ we call the interval $\makeset{i\in [1,\infty]}{\sigma_i=p}$ the \emph{$p$-interval} for $\sigma$.  
Let $P$ be any nonempty finite subset of $\P$. For each $p\geq 0$, define $M_p$ to be the poset consisting of all nonempty $p$-intervals for all $\sigma\in P$, ordered by inclusion. We define the \emph{$p$-profile} of $P$ to be the set 
$$\pr_p(P)=\makeset{(p,I)}{I \in \M_p\text{ is maximal} }.$$
Lastly, we define the \emph{profile} of $P$ to be the set
$$\pr(P):= \bigcup_{p\geq 0}\pr_p(P).$$
\end{definitions}

We pause to highlight some nuances in these definitions.   First, observe that a $p$-interval is infinite if and only if $p=0$.  Also,  a partition $\sigma$ has $p$-interval equal to $\emptyset$ if and only if $p$ is not a part in $\sigma$.  Lastly, if $p$ is not a part in any of the partitions in $P$, then $\pr_p(P) = \emptyset$.  

\begin{example}
If we take 
$$P = \{(4,2,2,2,1,1), (4,4,4,2,1,1,1), (4,4,2,2,2,2,1,1)\},$$
then 
\begin{align*}
\pr_4(P) &= \{(4,[1,3])\}\\
\pr_3(P) &= \emptyset\\
\pr_2(P) &= \{(2,[2,4]), (2,[3,6])\}\\
\pr_1(P) &= \{(1,[5,7]), (1,[7,8])\}\\ 
\pr_0(P) &= \{(0,[7,\infty] \}.
\end{align*}

\end{example}

\begin{definitions}
We say nonempty finite sets $P,Q\subseteq\P$ are \emph{profile-equivalent} and write $P\sim_{pr} Q$ if $P$ and $Q$ have the same profile.  In this way we obtain an equivalence relation on the set of nonempty finite subsets of $\P$. We denote the set of its equivalence classes by $\sP$. We call these classes \emph{profile classes}.  In the same way, we define an equivalence relation on the set of nonempty subsets of $\P(h,k)$, and denote the set of its equivalence classes by $\sP(h,k)$. 

\end{definitions}
Note that if $P\subseteq \P(h,k)$ then all elements of the profile class of $P$ are subsets of $\P(h,k)$, so the element of $\sP$ determined by $P$ coincides with the element of $\sP(h,k)$ determined by $P$.

The significance of profiles comes from the fact that 
\begin{equation}\label{eq:same profile same join}
\boxed{\quad P\sim_{pr} Q \Longleftrightarrow \vee_P = \vee_Q.\quad }
\end{equation}
This is established in the sequel by Theorem~\ref{thm:Closures, Profiles, and Maxima}.  In light of this result and our definition of a profile class, it makes sense to introduce the following notation.  
\begin{notation}
For any profile class $\C$, set
$$\vee_{\C}:=\vee_P\myand \Pr(\C):= \Pr(P),$$ 
for any $P\in \C$.	
\end{notation}

We are now in position to transform $F_{\mu,k}$ a bit further.   Starting with Equation~\ref{eq:F1}, we now have, for any $h,k >0$ and nonempty partition $\mu$ with $h\geq h_\mu$,
\begin{align}\label{eq:F2}
F_{\mu,k}(q)  &= \frac{1}{(1-q)\cdots (1-q^{k+w_\mu})}\sum_{\substack{P\subseteq \P(h,k)\\P\neq \emptyset}} (-1)^{|P|+1}  q^{|\vee_P(\mu)|}\nonumber\\
&= \frac{1}{(1-q)\cdots (1-q^{k+w_\mu})}\sum_{\C\in \sP(h,k)}\left(\sum_{P\in \C} (-1)^{|P|+1}  q^{|\vee_P(\mu)|}\right)\nonumber\\
&=\frac{1}{(1-q)\cdots (1-q^{k+w_\mu})} \sum_{\C\in \sP(h,k)} q^{|\vee_{\C}(\mu)|}\left(\sum_{P\in \C} (-1)^{|P|+1} \right) 	
\end{align}
where the last equality follows by the (yet unproven) fact given in (\ref{eq:same profile same join}).

Computational evidence suggests that the inner sum in (\ref{eq:F2}) is either $\pm 1$ or 0.  The profile classes for which the sum does not vanish can be nicely characterized as those having  a ``staircase" shape.  

\begin{definitions}
Let $h,k > 0$.  A set of the form
$$S=\{(p_1,[a_1,b_1]),(p_2,[a_2,b_2]),\ldots, (p_{s+1},[a_{s+1},b_{s+1}])\}$$
where
\begin{enumerate}[label=(\roman*),font=\itshape]
	\item $k=p_1>p_2>\cdots > p_{s+1} = 0$,
	\item $a_1=1$, $a_2\neq 1$,  $a_i\leq b_i$, and $b_{s+1} = \infty$,
	\item either $a_i = b_{i-1}$ or $a_i=b_{i-1}+1$
\end{enumerate}
is called an $(\ell,k)$-\emph{staircase} where $\ell = b_s$ is the \emph{length} of $S$. The set of all such staircases with length at most $h$ is denoted by $\stair(h,k)$.   An element $(p_i, [a_i,b_i])$ of $S$ is called \emph{left-overlapping} if $a_i = b_{i-1}$.

\end{definitions}

\begin{example}\label{ex:staircases}
Consider the set
$$P =\{(6,6,6,5,5,3,3),(6,6,5,5,5,3,2,2), (6,6,6,5,4,3,3)\}.$$
It has a staircase profile since 
$$\pr(P)=\{(6,[1,3]), (5,[3,5]), (4,[5,5]), (3,[6,7]), (2,[7,8]), (0,[8,\infty])\}.$$
Representing this graphically immediately suggests why the term ``staircase" is used:

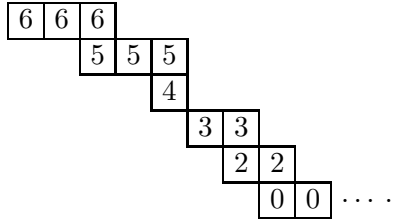
\begin{center}
\ytableausetup{boxsize=1.2em}
\begin{ytableau}
 6 & 6 & 6 \\
 \none & \none & 5 & 5 & 5\\
 \none & \none & \none & \none & 4  \\
  \none &  \none & \none & \none & \none & 3&3  \\
 \none & \none &\none &\none &\none &\none &2&2\\
 \none & \none &\none &\none &\none &\none &\none &0&0\\
\end{ytableau}\ \raisebox{-62pt}{$\ldots$\ .}
\end{center}

Before continuing we remark on our definition of length for a staircase.  The length of the above staircase is 8 which is precisely the largest height among all the partitions in $P$.  In fact, it is easily seen that this holds for any set of partitions with a staircase profile.  
\end{example}

Our next lemma states that all staircases $S\in\stair(h,k)$ are realized as profiles. Although we leave its proof to the reader, we note that our stipulation that the length of $S$ is no more than $h$ implies that  the elements of a set of partitions realizing $S$ must have height  at most $h$.  Furthermore, the fact that $p_1=k$ and $a_2\neq 1$ means that the largest part in each partition in such a set is $k$.  

\begin{lemma}\label{lem:all stair realized}
We have the following equality of sets: $$\stair(h,k)=\makeset{\pr(\C)}{\C\in \sP(h,k) \text{ has a staircase profile} }.$$
\end{lemma}

In light of Lemma~\ref{lem:all stair realized} it becomes natural to introduce another bit of notation.
\begin{notation}
For any $S\in\stair(h,k)$ we set
$\vee_{S}:=\vee_{\C}$
where $\C$ is the unique profile class whose profile is $S$.  	
\end{notation}

\begin{definitions}
	For any $\C\in \sP$ that has a staircase profile we say a set $M\subseteq \pr(\C)$  is an \emph{overlapping segment} provided the elements of $M$ can be ordered as
	$$(p_1, I_1), (p_2, I_2),(p_3, I_3),\ldots, (p_r, I_r),$$
	where $p_i>p_{i+1}$ and $I_{i+1}$ left-overlaps $I_i$.  We further define $\seg(\C)$ to be the number of maximal overlapping segments in $\pr(\C)$.  
\end{definitions}

The maximal overlapping segments in the example above are 
$$M_1 =\{ (6, [1,3]), (5, [3,5]), (4, [5,5])\} \myand M_2 = \{(3, [6,7]), (2, [7,8]), (0, [8,\infty])\}.$$
As in this example, we note that the maximal overlapping segments partition any staircase.

The relevance of staircases to the evolution of $F_{\mu,k}$ is given by Lemma~\ref{lem:only staircases survive} which establishes the following equality:
\begin{equation}\label{eq:outline:only staircases survive}
\boxed{
\quad\sum_{P\in \C}(-1)^{|P|} = 
	\begin{cases}
		(-1)^{|\pr(\C)|+\seg(\C)+1} &\text{ if $\C$  has a staircase profile}\\
		0 &\text{ otherwise}. 
	\end{cases}\quad}	
\end{equation}
Applying (\ref{eq:outline:only staircases survive}) along with Lemma~\ref{lem:all stair realized} to our formulation of  $F_{\mu,k}$ in Equation~\ref{eq:F2} we obtain, for any $h,k >0$ and nonempty partition $\mu$ with $h\geq h_\mu$,
\begin{align}\label{eq:F3}
	F_{\mu,k}(q) 	&=\frac{1}{(1-q)\cdots (1-q^{k+w_\mu})} \sum_{\C\in \sP(h,k)} q^{|\vee_{\C}(\mu)|}\left(\sum_{P\in \C} (-1)^{|P|+1} \right)\nonumber 	\\
	&=\frac{1}{(1-q)\cdots (1-q^{k+w_\mu})} \sum_{S\in \stair(h,k)} (-1)^{|S| +\seg(S)}q^{|\vee_{S}(\mu)|}.
\end{align}

Our final step is to translate staircases into what  we call \emph{augmented structures}.  This final step allows us to describe $F_{\mu,k}$ directly in terms of the Ferrers board of $\mu$ and the parameter $k$ without appealing to the function $\vee_S$ and (implicitly) to partitions that contain $\mu$.

\begin{definition}  Fix $h,k > 0$ and a nonempty partition $\mu$. We define a \emph{$(\mu,h,k)$-augmented structure} to be a 3-tuple of the form $(\mu,\lambda,\off)$ where $\lambda, \off\in \P$,
$$h_{\lambda} \leq h, \ \ h_{\off} \leq h + w_\lambda\ \ \textrm{and}\ \   w_{\lambda}+w_{\off}=k,$$
	and each column of $\lambda$ has length at least 2. Let $\aux(\mu,h,k)$ denote the set of $(\mu,h,k)$-augmented structures. We define the \emph{weight} of such a structure to be 
\begin{equation}\label{eq:weight of aux}
|(\mu, \lambda,\off)| := |\mu|+|\lambda|+|\off|+ \frac{a(a-1)}{2}+\sum_{i=1}^{a} \mu_{\lambda^*_i},
\end{equation}
where $a = w_\lambda$ and $\lambda^*_i$ denotes the length of the $i$th column of $\lambda$.  In this setting, the partition $\lambda$ is called an $L$-\emph{partition} and $\off$ is called an \emph{offset partition}.
\end{definition}

Note that one of $\lambda, \off$ may be empty in an augmented structure.

To shed light on our formula for weight and our choice of terminology, consider the following example where we take $h=7$:
$$\mu = (4,3,3,2,2,1),\quad \lambda=(4,4,3,3,2,1,1),\quad\off = (4,2,2,1).$$
Draw the Ferrers boards for $\mu$, $\lambda$, and $\off$ so that the board for $\mu$ is right justified as shown:
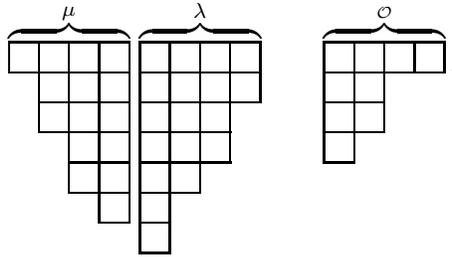
\ytableausetup{boxsize=1em}
$$\overbrace{\ydiagram{4,1+3,1+3,2+2,2+2,3+1}}^{\mu}\, \overbrace{\ydiagram{4,4,3,3,2,1,1}}^{\lambda}\qquad \overbrace{\ydiagram{4,2,2,1}}^{\off}.$$
Now consider all (reversed) $L$'s formed by first choosing the boxes in some column of $\lambda$ and then choosing all boxes to the left of the bottom box of this column including those in $\mu$.   For example, if we highlight the $L$'s anchored in the 1st and 3rd columns we have
$$\ydiagram{4,1+3,1+3,2+2,2+2,3+1}*[*(red)]{0,0,0,2+2}\, \ydiagram{1+3,1+3,1+2,1+2,1+1}*[*(purple)]{0,0,0,1}*[*(red)]{2+1,2+1,2+1,1+2}*[*(blue)]{1,1,1,1,1,1,1}\qquad \ydiagram{4,2,2,1}.$$
The number of boxes in the red $L$ is $\lambda^*_3 + (3-1) + \mu_{\lambda^*_3} =  8$ and the number of boxes in the blue $L$ is $\lambda^*_1 + (1-1) + \mu_{\lambda^*_1} =  7$.

In light of this, we interpret the second, fourth, and fifth  terms in (\ref{eq:weight of aux}) as giving the total number of boxes involved in all possible such $L$'s.  
Our main result relating staircases with augmented structures (see Lemma~\ref{lem:bijection between stair and aux}) is the following statement.

\medskip
\begin{equation}\label{eq:staircases to aug}
\boxed{
\begin{minipage}[c]{.9\textwidth}
Fix $h,k > 0$ and a partition $\mu$ so that $h\geq h_\mu$. Then there exists a bijection between $\stair(h,k)$ and $\aux(\mu,h,k)$ so that if $S\in \stair(h,k)$ maps to $(\mu,\lambda,\off)\in\aux(\mu,h,k)$, then
	\begin{enumerate}[label=(\roman*),font=\itshape]
		\item $w_{\lambda} = |S|-\seg(S)$, and
		\item $|\join_S (\mu)| = |(\mu,\lambda,\off)|$.
	\end{enumerate} 
\end{minipage}}
\medskip
\end{equation}

We pause to point out a subtlety in the definition of augmented structures.  In our definition we insist that every column of $\lambda$ has length at least 2.  This is required in order for staircases to correspond with our augmented structures.  In particular, the requirement that $a_2\neq 1$ in the definition of staircases translates into this condition on augmented structures.  See Lemma~\ref{lem:bijection between stair and aux} for details.

Using Lemma~\ref{lem:bijection between stair and aux} we arrive at our final description of $F_{\mu,k}$.

\begin{theorem}\label{thm:GF of F_{mu,k}}
Fix $h,k >0$ and a nonempty partition $\mu$ so that $h\geq h_\mu$.  Then
$$F_{\mu,k}(q)= \frac{1}{(1-q)\cdots (1-q^{k+w_\mu})} \sum_{(\mu,\lambda,\off)\in \aux(\mu,h,k)} (-1)^{w_{\lambda}}q^{ |(\mu,\lambda,\off)|}.$$
\end{theorem}

\begin{proof}
Applying Lemma~\ref{lem:bijection between stair and aux} to Equation~\ref{eq:F3} and observing that $|S|-\seg(S)$ has the same parity as $|S|+\seg(S)$ we obtain
\begin{align*}
	F_{\mu,k}(q)&=\frac{1}{(1-q)\cdots (1-q^{k+w_\mu})} \sum_{S\in \stair(h,k)} (-1)^{|S| +\seg(S)}q^{|\vee_{S}(\mu)|}\\
	&= \frac{1}{(1-q)\cdots (1-q^{k+w_\mu})} \sum_{(\mu,\lambda,\off)\in \aux(\mu,h,k)} (-1)^{w_{\lambda}}q^{|(\mu, \lambda,\off)|}.
\end{align*} 	
This proves the theorem.   
\end{proof}

  \section{Rook equivalence and Wilf equivalence}\label{sec:equivalence}
  
  Using Theorem~\ref{thm:GF of F_{mu,k}} we can now establish connections between the concepts of rook equivalence and Wilf equivalence of integer partitions.  We begin by establishing necessary and sufficient conditions for width-Wilf equivalence. It is clear that if two partitions are width-Wilf equivalent then they must have the same weight.
  
  \begin{theorem}\label{thm:width Wilf}  For any nonempty partitions $\mu$ and $\tau$ that have the same weight and width, the following are equivalent:
  
  \begin{enumerate}[label=(\roman*)]
      \item $\mu$ and $\tau$ are width-Wilf equivalent
      \item $F_{\mu,1}=F_{\tau,1}$
      \item $\mu$ and $\tau$ are rook equivalent.
  \end{enumerate}
  
  \end{theorem}
  
   \begin{proof} We show that $(i)\Rightarrow (ii)\Rightarrow (iii)\Rightarrow (i)$.
   
   If $\mu$ and $\tau$ are width-Wilf equivalent and  $w_{\mu}=w_{\tau}=w$, then $\mu$ and $\tau$ have the same number of extensions of width $w+1$ and any specified weight, so $F_{\mu,1}=F_{\tau,1}$.  This proves that $(i)\Rightarrow (ii)$.
   
   To prove that $(ii)\Rightarrow (iii)$, assume that $(ii)$ holds, and let $h=\max\{h_{\mu},h_{\tau}\}$. Then Theorem~\ref{thm:GF of F_{mu,k}} implies that
   
   $$\sum_{(\mu,\lambda,\off)\in \aux(\mu,h,1)} (-1)^{w_{\lambda}}q^{ |(\mu,\lambda,\off)|}=\sum_{(\tau,\lambda,\off)\in \aux(\tau,h,1)} (-1)^{w_{\lambda}}q^{ |(\tau,\lambda,\off)|}.$$
   
   By abuse of notation, denote the augmented structures in $\aux(\mu,h,1)$ by $\lambda_2,\ldots, \lambda_{h}$ and $\off_1,\ldots, \off_{h}$, where in $\lambda_i$ the offset partition is empty and the $L$-partition has a single column of weight $i$, and in $\off_j$ the $L$-partition is empty and the offset partition has a single column of weight $j$. If $\mu$ and $\tau$ have weight $n$, then
       $|\off_j|= n+j$ for each $j$, and $|\lambda_i|=n+\mu_i+i$. Therefore
       
        $$-\sum_{2\leq i\leq h}q^{n+\mu_i+i}+ \sum_{1\leq j\leq h} q^{n+j}=
      -\sum_{2\leq i\leq h}q^{n+\tau_i+i}+
      \sum_{1\leq j\leq h}q^{n+j},$$
       and it follows that the multisets 
       
       $$\{\mu_2+2,\ldots, \mu_h+h\}\myand \{\tau_2+2,\ldots, \tau_h+h\}$$ are equal.  Since $\mu_1=\tau_1$ by assumption, the multisets
       
        $$\{\mu_1+1,\ldots, \mu_h+h\} \myand  \{\tau_1+1,\ldots, \tau_h+h\}$$
        are equal, and this condition is equivalent to the rook equivalence of $\mu$ and $\tau$ (see~\cite{loehr2011bijective}, Theorem 12.10, which was proved in  \cite{foata1970rook} in a slightly different form). This concludes the proof that $(ii)\Rightarrow (iii)$.
        
        To prove that $(iii)\Rightarrow (i)$, assume that $(iii)$ holds, and again let $h=\max\{h_{\mu},h_{\tau}\}$ and let $\mu$ and $\tau$ have weight $n$. It is clear that $\mu$ and $\tau$ have the same number of extensions of width $w$ and any specified weight, obtained by adding the same rows of length at most $w$ to $\mu$ and $\tau$.  So to prove that $\mu$ and $\tau$ are width-Wilf equivalent it will suffice to show that  $F_{\mu,k}=F_{\tau,k}$ for all $k\geq 1$.  Since $\mu$ and $\tau$ have the same width, what we need to show is that

   $$\sum_{(\mu,\lambda,\off)\in \aux(\mu,h,k)} (-1)^{w_{\lambda}}q^{ |(\mu,\lambda,\off)|}=\sum_{(\tau,\lambda,\off)\in \aux(\tau,h,k)} (-1)^{w_{\lambda}}q^{ |(\tau,\lambda,\off)|}$$ for each $k\geq 1$.  To show this it suffices to establish a bijection from $\aux(\mu,h,k)$ onto $\aux(\tau,h,k)$ that preserves the weight of each augmented structure and the width of its $L$-partition.
   
   Because $\mu$ and $\tau$ are rook equivalent and have the same width, the multisets
    
     $$\{\mu_2+2,\ldots, \mu_h+h\} \ \ \textrm{and}\ \  \{\tau_2+2,\ldots, \tau_h+h\}$$ are equal, by the result cited above. Therefore there exists a bijection $\varphi:[2,h]\rightarrow [2,h]$ such that $\mu_i+i=\tau_{\varphi(i)}+\varphi(i)$ for all $i\in [2,h]$.  We obtain the desired bijection from $\aux(\mu,h,k)$ onto $\aux(\tau,h,k)$ by sending each $(\mu,\lambda,\off)$ to $(\tau,\lambda',\off)$, where the lengths of the columns of $\lambda'$ are obtained by applying $\varphi$ to the lengths of the columns of $\lambda$.

    \end{proof}
    
    If partitions $\mu$ and $\tau$ are rook equivalent but have different widths, then although $\mu$ and $\tau$ cannot be width-Wilf equivalent, they  must be Wilf equivalent.  To prove this, we will use a criterion for rook equivalence different from the multiset criterion employed above.
    
    Let $F$ be a Ferrers board, and let $(i,j)$ denote the square in the $i$th row from the top and the $j$th column from the left. If $(i,j)$ has the property that replacing the subboard $\makeset{(x,y)}{x\geq i, \ y\geq j}$ by its conjugate yields a Ferrers board, then, following Foata and Sch\"{u}tzenberger \cite{foata1970rook} we call this new board  the \emph{(i,j)-transform} of $F$.  We will denote it by $F_{i,j}$. Foata and Sch\"{u}tzenberger proved (\cite{foata1970rook}, Corollary 6, Lemma 9 and Theorem 11) that boards $F_1$ and $F_2$ are rook equivalent if and only if one can be obtained from the other by a series of these transformations.
    
    \begin{theorem}\label{thm:rook and Wilf}  If partitions $\mu$ and $\tau$ are rook equivalent, then they are Wilf equivalent.
    \end{theorem}
    
    \begin{proof}
    
    We view $\mu$ and $\tau$ as Ferrers boards.  If $\mu$ and $\tau$ are rook equivalent, then by the criterion of Foata and Sch\"{u}tzenberger we can obtain $\tau$ from $\mu$ by a series of $(i,j)$-transformations. It will suffice to show that, for each $(i,j)$-transformation, $F$ and $F_{i,j}$ are Wilf equivalent.
    
    \medskip
    Case 1: $i>1$.  In this case, $F$ and $F_{i,j}$ are rook equivalent and have the same width, so by Theorem~\ref{thm:width Wilf} they are width-Wilf equivalent and therefore Wilf equivalent.
    
    \medskip
    Case 2: $(i,j)=(1,1)$. In this case every extension of $F$ corresponds (via a conjugation) to an extension of $F_{1,1}$ of the same weight.
    
    \medskip
    Case 3: $i=1$ and $j>1$. In this case $F_{i,j}$ can be obtained in three steps, by taking $F_{1,1}$, then taking the $(j,i)$-transform of the result, and then taking the $(1,1)$-transform again.  Each step yields a result Wilf equivalent to $F$, as we see by arguing as in Case 2 for the first and last steps, and as in Case 1 for the middle step (using $j>1$).
       
     \end{proof} 
     
The converse of Theorem~\ref{thm:rook and Wilf} is also true.  We prove it in \cite{Bloom:On-cr2017}, but we remark here that, for the proof,
 it suffices to show that two Wilf equivalent partitions with distinct parts must be identical. For if any partitions $\mu$ and $\tau$ are Wilf equivalent, then (\cite{foata1970rook}, Theorem 11) $\mu$ and $\tau$ are rook equivalent to partitions $\mu'$ and $\tau'$ with distinct parts. By Theorem~\ref{thm:rook and Wilf}, $\mu'$ and $\tau'$ are Wilf equivalent, so, assuming that this forces them to be equal, we see that $\mu$ and $\tau$ are rook equivalent.

\section{Profiles}\label{sec:profiles}
In this section we further explore the idea of profiles, first introduced in Section~\ref{sec:GF for P(mu,k)}.  This section also contains the proofs of (\ref{eq:same profile same join}) and  (\ref{eq:outline:only staircases survive}). 

In Section~\ref{sec:GF for P(mu,k)}, we claimed (without proof), that two nonempty finite sets of partitions $P,Q\subseteq \P$ have the same profile if and only if the mappings $\vee_P,\vee_Q:\P\to \P$ are equal.  Proving this is our first order of business.  We start with a couple of definitions.  

\begin{definitions}
	Suppose $\alpha, \beta \in \P$.  Provided $\alpha_i > \beta_{i+1}$, we define
	$$\alpha\star_i \beta = (\alpha_1,\ldots, \alpha_i,\beta_{i+1},\ldots) \in \P.$$
We say that $\alpha\star_i \beta$ is obtained from $\alpha$ and $\beta$ by \emph{splicing}.  

For any set $P\subseteq \P$ we denote by $\cl(P)$ the closure of $P$ under splicing (using all $i>0$). 
 
\end{definitions}

\begin{example}\label{ex:closure}  Let $P = \{ (2,1,1), (2,2,1,1) \}$ and note that 
$$\cl(P) =\{(2,1,1), (2,2,1),(2,2,1,1)\}.$$
Comparing the profiles of $P$ and $\cl(P)$, we find that 
$$\pr(P) = \{(2,[1,2]),(1,[2,3]),(1,[3,4]), (0,[4,\infty])\}=\pr(\cl(P)).$$
This suggests that the profile of a finite set is invariant under splicing.   The next two results prove this and much more.  Before launching into these proofs, we remark on a subtle point in the definition of splicing.  In the definition, $\alpha\star_i \beta$ is defined if and only if $\alpha_i>\beta_{i+1}$.  Of course the operation, a priori, could have been defined with $>$ replaced by $\geq$.  This would certainly yield a partition, but it would not preserve the profile of a set.  For example, under this alternative definition the ``closure" of $P$ would be the set 
$$\{(2,1,1), (2,2,1),(2,2,1,1), (2,1,1,1)\},$$
whose profile is $\{(2,[1,2]),(1,[2,4]), (0,[4,\infty])\}$.
\end{example} 

We now turn our attention to proving (\ref{eq:same profile same join}) by first proving a technical lemma that is of interest in its own right.

\begin{lemma}\label{lem:adding a partition}
    If $P \subseteq \P$ is finite and nonempty  and $\gamma\in \P$ then the following are equivalent:
\begin{enumerate}[label=(\roman*)]
\item $\gamma\in \cl(P)$
\item $\join_P= \join_{P\cup \{\gamma\}}$
\item $\pr(P)=\pr(P\cup \{\gamma\})$.
\end{enumerate}

\end{lemma}

\begin{proof}
   We show that $(i)\Rightarrow (ii)\Rightarrow(iii)\Rightarrow (i)$.
   
   To begin, we first show that if $\mu, \alpha, \beta\in \P$, and $\gamma=\alpha \star_{\ell} \beta$, then
   $$(\mu+\alpha) \join  (\mu+\beta)= (\mu+\alpha) \join  (\mu+\beta) \join  (\mu+\gamma).$$
   To see this, we first note that the first $\ell$ parts of $\mu+\gamma$ are parts of $\mu+\alpha$, and the remaining parts of $\mu+\gamma$ are parts of $\mu+\beta$, so $\mu+\gamma$ does not introduce new parts into $(\mu+\alpha)\vee (\mu+ \beta)$.
Now suppose for a contradiction that $\mu+\gamma$ has some part $r$ that has a higher multiplicity in $\mu+\gamma$ than it does in either $\mu+\alpha$ or $\mu+\beta$. Then $\makeset{i}{(\mu+\gamma)_i=r}$ must include both $\ell$ and $\ell +1$, so $(\mu+\gamma)_{\ell}=(\mu+\gamma)_{\ell +1}$, a contradiction since $\alpha_{\ell}> \beta_{\ell+1}$.  

Now to prove that $(i)\Rightarrow (ii)$, suppose $\gamma\in\cl(P)$.
Then there exists a sequence of partitions $\delta_1,\ldots,\delta_t$ such that each $\delta_j$ is obtained by applying the splicing operation to elements of $P$ and the previous $\delta_i$'s, and $\delta_t=\gamma$.  Using the result of the preceding paragraph repeatedly, we have
$$\join_P(\mu)=\join_{P\cup \{\delta_1\}} (\mu)=\cdots = \join_{P\cup \{\delta_1,\ldots, \delta_{t}\}} (\mu),$$
so $\join_P (\mu)=\join_{P\cup \{\gamma\}} (\mu)$.  (If $\vee_P(\mu)$ does not change when we add all of the $\delta_i$'s to $P$, then it does not change when we only add $\delta_t=\gamma$ to $P$.)

To prove that $(ii)\Rightarrow (iii)$, assume $\join_P (\mu)= \join_{P\cup \{\gamma\}} (\mu)$ for every $\mu\in\P$. We show that every interval in $\pr(\{\gamma\})$ is contained in some interval in $\pr(P)$, so $\pr(P)=\pr(P\cup \{\gamma\})$. 
To see this, suppose $r$ is some part of $\gamma$ and let $I=[a,b]$ be its $r$-interval.  So 
$$\gamma = (\gamma_1,\ldots, \gamma_{a-1}, \underbrace{r,r,\ldots, r}_{|I|}, \gamma_{b+1},\ldots).$$
Next set $\displaystyle k = \max_{\sigma\in P\cup\{\gamma\}} w_\sigma$ and define the partition $\mu$ to be 
$$\mu = (\underbrace{4k,4k,\ldots, 4k}_{a-1},\underbrace{2k,2k,\ldots 2k}_{|I|}).$$
Then $(\mu+\gamma)_i=2k+r$ for all $i\in I$. Because $\vee_P(\mu) = \vee_{P\cup \{\gamma\}} (\mu)$, there must exist some $\sigma\in P$ such that the multiplicity of $2k+r$ in $\mu+\sigma$ is at least $|I|$. Since $r\leq k$ and each part of $\sigma$ is at most $k$, $(\mu+\sigma)_i$ can be $2k+r$ only when $i\in I$, so $(\mu+\sigma)_i=2k+r$ for all $i\in I$.  Thus $I\subseteq \makeset{i}{\sigma_i=r}$, and this implies that $I$ is contained in some interval in $\pr(P)$.

To prove that $(iii)\Rightarrow (i)$, suppose $\pr(P)=\pr(P\cup \{\gamma\})$. Denote the elements of $\pr(\gamma)$ by
$$(p_1,I_1),\ldots, (p_m, I_m)$$
so that $k=p_1>p_2>\cdots >p_m =0$.  Let $b_i$  be the right endpoint of $I_i$.  (Note that in this case $b_i+1$ is the left endpoint of $I_{i+1}$ for all $i<m$.)
By our assumption that $\pr(P) = \pr(P\cup \{\gamma\})$, it follows that for each $1\leq i\leq m$ there exists some $\sigma^{(i)}\in P$ whose $p_i$-interval contains $I_i$.  Now define
$$\gamma' = \sigma^{(1)}\star_{b_1} \sigma^{(2)} \star_{b_2} \cdots \star_{b_{m-1}} \sigma^{(m)}$$
where we evaluate this expression from left to right.  We have $\pr(\{\gamma'\}) = \pr(\{\gamma\})$, so $\gamma = \gamma'\in \cl(P)$.  
\end{proof}

We are now ready to prove (\ref{eq:same profile same join}).   

\begin{theorem}\label{thm:Closures, Profiles, and Maxima}
  If $P,Q\subseteq \P$ are nonempty finite sets then the following are equivalent:
\begin{enumerate}[label=(\roman*)]
\item $\cl(P)=\cl(Q)$
\item $\pr(P)=\pr(Q)$
\item for every $\mu\in \P$, we have $\vee_{P}(\mu) = \vee_{Q}(\mu)$.
\end{enumerate}
\end{theorem}

\begin{proof}
We show that $(i)\Leftrightarrow (ii)$ and $(i)\Leftrightarrow (iii)$.

To show that $(i)\Rightarrow (ii)$, let the elements of $\cl(P)$ be $\beta_1,\ldots,\beta_{\ell}$.  Then by repeated use of the preceding lemma we have

$$\pr(P)=\pr(P\cup\{\beta_1\})=\cdots =\pr(P\cup\{\beta_1,\ldots,\beta_{\ell}\})=\pr(\cl(P)),$$
and likewise $\pr(Q)=\pr(\cl(Q))$.  Thus if $\cl(P)=\cl(Q)$ then $\pr(P)=\pr(Q)$.

To show that $(ii)\Rightarrow (i)$, assume $\pr(P)=\pr(Q)$.  Then $\pr(P\cup \{\gamma\})= \pr(Q\cup \{\gamma\})$ for every $\gamma\in \P$.  Thus, using the preceding lemma, we have 
$$\gamma\in \cl(P)\Leftrightarrow \pr(P)=\pr(P\cup \{\gamma\})\Leftrightarrow \pr(Q)=\pr(Q\cup \{\gamma\})\Leftrightarrow \gamma\in \cl(Q).$$  
So $\cl(P)=\cl(Q)$.

A proof that $(i)\Leftrightarrow (iii)$ can be obtained from our proof that $(i)\Leftrightarrow (ii)$ by replacing $\pr(X)$ by $\join_X (\mu)$ for all subsets $X$ of $\P$ involved.

\end{proof}

In light of this theorem it is natural to introduce the following notation.
\begin{notation}
Let $\C\in \sP$, so that for all $P,Q\in \C$ we have $\pr(P) = \pr(Q)$. Define
$\cl(\C):= \cl(P)$  for any $P\in\C$.
\end{notation}

\begin{remark}
For any $P\in \C$ we have $\cl(\cl(P))=\cl(P)$, so by our theorem, $\pr(\cl(P))=\pr(P)$ and thus $\cl(\C)\in \C$. We also note that 
$$\cl(\C) \subseteq \bigcup_{P\in \C} P \subseteq  \bigcup_{P\in \C} \cl(P) = \cl(\C),$$
where the first inclusion follows since 
$\cl(\C)) \in \C $ and the equality on the right follows since $\cl(\C)= \cl(P)$ for every $P\in \C$.  This gives us another description of the set $\cl(\C)$, as the union of all members of $\C$. 
\end{remark}

Next, we turn our attention to the proof of (\ref{eq:outline:only staircases survive}), which appears in Lemma~\ref{lem:only staircases survive}. Some preliminary results are required before tackling that lemma.

The next definition allows us to frame the proof of Lemma~\ref{lem:only staircases survive} in a convenient graph-theoretic context.

\begin{definition}
	Given a profile class $\C$ let $G(\C)$ be the bipartite graph whose partite sets are the set of partitions $\cl(\C)$ and the profile $\pr(\C)$,  such that $\alpha\in\cl(\C)$ is adjacent to $(p,I)\in\pr(\C)$ if and only if $I$ is the $p$-interval in $\alpha$.  
	
	In keeping with traditional graph theory notation, we denote the set of neighbors of a vertex $w$ in $G(\C)$ by $N(w)$ and define 
 $$N(W) = \bigcup_{w\in W} N(w)$$
 for any collection of vertices $W$ in $G(\C)$.  
\end{definition}

\begin{example}
In the case that our profile class is as in Example~\ref{ex:closure} the corresponding bipartite graph is: 
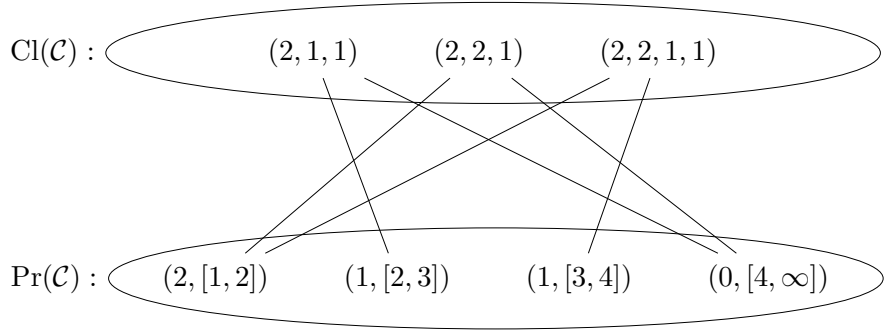
\begin{center}

	\begin{tikzpicture}[thick,
  every fit/.style={line width = 0.1mm, ellipse,draw,inner sep=4pt,text width=7cm}]

% the vertices of U
\begin{scope}[start chain=going right,node distance=7mm]
	\node[on chain] (a)  {$(2,[1,2])$};
	\node[on chain] (b)  {$(1,[2,3])$};
	\node[on chain] (c)  {$(1,[3,4])$};
	\node[on chain] (d)  {$(0,[4,\infty])$};
\end{scope}

% the vertices of V
\begin{scope}[xshift=13mm,yshift=3cm,start chain=going right,node distance=7mm]
\node[on chain] (1) {$(2,1,1)$};
\node[on chain] (2) {$(2,2,1)$};
\node[on chain] (3) {$(2,2,1,1)$};

\end{scope}

% the set U
\node [fit=(1) (3),label=left:$\cl(\C):$] {};
% the set V
\node [fit=(a) (d),label=left:$\pr(\C):$] {};

% the edges
\draw[line width = 0.1mm] (1) -- (b);
\draw[line width = 0.1mm] (1) -- (d);
\draw[line width = 0.1mm] (2) -- (a);
\draw[line width = 0.1mm] (2) -- (d);
\draw[line width = 0.1mm] (3) -- (a);
\draw[line width = 0.1mm] (3) -- (c);
\end{tikzpicture}
\end{center}	
\end{example}

In stating the next lemma, we adopt the convention that a sum over the empty set is $0$.

\begin{lemma}\label{lem:bipartite graph1}
Let $G$ be a finite bipartite graph with (nonempty) partite sets $A$ and $B$. Then we have
	$$\sum_{\substack{S\subseteq A\\ N(S) = B}} (-1)^{|S|} = \sum_{\substack{S\subseteq B\\ N(S) = A}} (-1)^{|S|}.$$
\end{lemma}
\begin{proof}
	The proof of this follows by a straightforward application of the Inclusion--Exclusion principle and the fact that 
\begin{align}\label{eq:power set}
\sum_{S\subseteq U} (-1)^{|S|} = \begin{cases}
1 &\text{ if } U = \emptyset\\
0 &\text{ if } U \neq \emptyset
 \end{cases}	
\end{align}
where $U$ is an arbitrary finite set.  Using these two tools we obtain the calculation
\allowdisplaybreaks
	\begin{align*}
		\sum_{\substack{S\subseteq A\\ N(S) = B}} (-1)^{|S|} &=\sum_{S\subseteq A} (-1)^{|S|} - \sum_{\substack{T\subseteq B\\T \neq \emptyset}}\sum_{\substack{S\subseteq A\\ N(S) \subseteq B\setminus T}} (-1)^{|S|+|T|+1}\\
		&=-\sum_{\substack{T\subseteq B\\T \neq \emptyset}}\sum_{\substack{S\subseteq A\\ N(S) \subseteq B\setminus T}} (-1)^{|S|+|T|+1}\\
		&=\sum_{\substack{T\subseteq B\\T \neq \emptyset}}\sum_{\substack{S\subseteq A\\ N(S) \subseteq B\setminus T}} (-1)^{|S|+|T|}\\
		&=\sum_{\substack{T\subseteq B\\T \neq \emptyset}}\sum_{\substack{S\subseteq A \setminus N(T)}} (-1)^{|S|+|T|}\\
				&=\sum_{\substack{T\subseteq B\\  N(T)=A}} (-1)^{|T|},
	\end{align*}
where the second equality follows since $A\neq \emptyset$.  This concludes our proof.
\end{proof}

\begin{definition}
Let $G$ be a finite bipartite graph with (nonempty) partite sets $A$ and $B$. Fix $X\subseteq A$.  We say a set $Y\subseteq B$ is $X$-\emph{minimal} provided that $X \subseteq N(Y)$ yet the neighbors of any proper subset of $Y$ do not contain all of $X$. 
\end{definition}

\begin{lemma}\label{lem:bipartite graph2}
	Let $G$ be a finite bipartite graph with (nonempty) partite sets $A$ and $B$.  Assume there exist distinct vertices $u,v\in B$ such that $N(u)\cap N(v)=\emptyset$ and no $A$-minimal set contains both $u$ and $v$. Then 
	$$\sum_{\substack{S\subseteq B\\ N(S) = A}} (-1)^{|S|} = 0.$$
\end{lemma}
\begin{proof}
		Define $$\sB = \makeset{S\subseteq B}{N(S)=A}.$$ To prove the lemma, it suffices to define an involution $f$ on the collection $\sB$ with the property that $S$ and $f(S)$ have opposite parity.  To help define our involution we say a set $S\in \sB$ is $u$-\emph{critical} provided that $N(S\setminus u) \neq A$.  (Note that a $u$-critical set $S$ must contain $u$.)  Now let
		
		$$f(S) = \begin{cases}
		S\cup \{u\} & \text{if } u\notin S\\
		S\setminus \{u\} & \text{if } u\in S \text{ and $S$ is not $u$-critical}\\
		S\cup \{v\} & \text{if } v\notin S \text{ and } $S$ \text{ is $u$-critical}\\
		S\setminus \{v\} & \text{if } v\in S \text{ and } $S$ \text{ is $u$-critical}.\\
	\end{cases}$$
	 
We must check that $f$ is well defined and an involution.  To begin, observe that if $S\in \sB$ falls in the first or third case, then clearly $f(S)\in \sB$.  If $S$ falls in the second case, then $f(S)\in \sB$ since $S$ is assumed not to be $u$-critical.  If $S$ falls in the fourth case, then $N(S\setminus \{v\})$ must equal $A$ for otherwise every suset $T\subseteq S$ such that $N(T)=A$ would have to contain $u$ and $v$, so we would have an $A$-minimal set that contains both $u$ and $v$.  Hence $f(S)\in \sB$ in this case as well. The four cases in the definition are exhaustive and mutually exclusive, by the parenthetical note preceding the definition of $f$. 

We now turn our attention to showing that $f$ is indeed an involution.  First observe that $f$ certainly alternates between sets $S$ that fall in the first case and those that fall in the second case.  If $S$ falls in the fourth case, then removing $v$ from the $u$-critical set $S$ clearly results in another $u$-critical set $f(S)$ that does not contain $v$.  On the other hand, if $S$ falls in the third case we must show that $f(S) = S\cup \{v\}$ is also $u$-critical.  The only way this might not happen is if $N(u)\cap N(v) \neq \emptyset$.  As this is not the case, we can conclude $f$ alternates between sets $S$ that fall in the third case and those that fall in the fourth case. 
\end{proof}

 \begin{lemma}\label{lem:vanishing profiles}
	Fix a profile class $\C$ and distinct profile elements $(p,I),(q,J)\in \pr(\C)$, where $I=[a,b]$, $J=[c,d]$, $p\geq q$, and $c< b$. Then 
	$$\sum_{P\in \C} (-1)^{|P|}=0.$$
 \end{lemma}

\begin{proof}
First observe that 
$$\sum_{P\in \C} (-1)^{|P|} = \sum_{\substack{P\subseteq \cl(\C)\\ N(P) = \pr(\C)}} (-1)^{|P|}=\sum_{\substack{X\subseteq \pr(\C)\\ N(X) = \cl(\C)}} (-1)^{|X|},$$
where the second equality comes from Lemma~\ref{lem:bipartite graph1} and the first follows from the fact that for all $P\subseteq \cl(\C)$, we have
$N(P) = \pr(\C)\Leftrightarrow \pr(P) = \pr(\C) \Leftrightarrow P\in \C$.  
The remainder of our proof is devoted to showing that the righthand side is $0$. We do this by showing that our graph $G(\C)$ contains two vertices $u,v\in\pr(\C)$ that satisfy the hypotheses of Lemma~\ref{lem:bipartite graph2}.  

To begin we choose
$$u = (p,[a,b])\myand v = (q,[c,d]),$$
so that $p\geq q$ and $c<b$ and $u\neq v$. A straightforward check verifies that $N(u)\cap N(v) =\emptyset$.  We leave the details to the reader.

Next we show that no $\cl(\C)$-minimal set contains both $u$ and $v$.  For a contradiction, assume that $M\subseteq \pr(\C)$ is a $\cl(\C)$-minimal set that contains both $u$ and $v$.  This means that 
$$ \underbrace{N(u) \setminus N(M-u)}_{X} \neq \emptyset \myand  \underbrace{N(v) \setminus N(M-v)}_{Y}\neq \emptyset,$$
where $M-x := M \setminus \{x\}$.  Now fix $\alpha\in X$ and $\beta \in Y$.  Clearly, $\alpha\neq \beta$ as $N(u)$ and $N(v)$ are disjoint.

We claim that there exists $\gamma\in \cl(\C)$ such that $\gamma\notin N(u)\cup N(v)$. To see this we consider the following two cases. Only these two are necessary since when $p=q$, the fact that neither $I$ nor $J$ is a subset of the other implies that, by switching the roles of $I$ and $J$, if need be, we may take $d<b$.  

\medskip
Case 1: $p>q$ and $d\geq b$
\medskip

In this case, define $\gamma = \alpha \star_{b-1} \beta\in \cl(\C)$ where splicing is permitted since $\beta_b = q$ and $\alpha_{b-1}\geq \alpha_b= p>q$.  In this case, $\gamma_b=\beta_b=q\neq p$, so $\gamma \notin N(u)$.  Similarly, $\gamma_c=\alpha_c\geq p> q$, so $\gamma\notin N(v)$.    

\medskip
Case 2: $p\geq q$ and $d< b$
\medskip

In this case, define $\gamma = \alpha \star_{d} \beta\in \cl(\C)$ where splicing is permitted since $\beta_{d+1} < q\leq p \leq \alpha_d$ because $\alpha_b=p$ and $d<b$.  In this case, $\gamma_b=\beta_b< q,$ so $\gamma_b< p$ and therefore $\gamma \notin N(u)$. Also, $\gamma_d=\alpha_d\geq p$ since $d< b$.  If $p> q$ then $\gamma_d\neq q$ so $\gamma\notin N(v)$. If $p=q$ then $c< a$ by definition of a profile.  Then $\gamma_c=\alpha_c> p=q$, so $\gamma\notin N(v)$.

In either case we have $\gamma \in \cl(\C) = N(M)$ with $\gamma \notin N(u)\cup N(v)$.  Since $M$ is $\cl(\C)$-minimal we must then have some $w\in M\setminus \{u,v\}$ with $\gamma \in N(w)$.  Taking either definition of $\gamma$, it follows that either $\alpha \in N(w)$ or $\beta\in N(w)$.  But the first option contradicts the fact that $\alpha\in X$  and the second contradicts the fact that $\beta \in Y$.  We conclude that no $\cl(\C)$-minimal set can contain both $u$ and $v$. 
\end{proof}

We turn our attention to staircases.  The reader might find it helpful to revisit Example~\ref{ex:staircases} in Section~\ref{sec:GF for P(mu,k)} before continuing. 

Our first task is to introduce an important collection of partitions associated with staircases.  An example of the construction described in this lemma is given immediately after its proof.  

\begin{lemma}
Fix some $\C\in \sP$ whose profile is a staircase.  Denote the maximal overlapping segments in this staircase as $M_1,\ldots, M_\ell$ and choose $(p_i,I_i)\in M_i$ for each $i\leq \ell$.  Then there is a partition $\sigma(\C, p_1,\ldots, p_\ell)\in \cl(\C)$ whose neighbors in $G(\C)$ are precisely  
$$(p_1,I_1),\ldots, (p_\ell,I_\ell).$$
\end{lemma}

\begin{proof}
	First, for any maximal overlapping segment $M$ with $(p,I)\in M$ we define
	$$\pi(M,p) = \bigcup_{\substack{(q,J)\in M\\ q>p\\|J|>1}} \{(q,J^-)\} \cup \{(p,I)\} \cup \bigcup_{\substack{(q,J)\in M\\ p>q\\ |J|>1}}\{(q,{^-J})\}$$
	where for any interval $I$ we define $I^-$ (respectively, $^-I$) to be the interval obtained by deleting its right (respectively, left) endpoint. Observe that $\pi(M,p)\cap M = \{(p,I)\}$. 
 
 Consider the set
$$\pi(M_1,p_1)\cup \cdots \cup \pi(M_\ell,p_\ell).$$
Observe that this set is the profile of some (unique) partition $\sigma$.  Furthermore, since $\pr(\cl(\C)\cup \{\sigma\})=\pr(\cl(\C))$, it follows from Lemma~\ref{lem:adding a partition} that $\sigma\in \cl(\C)$.

Our first observation now implies that
$$\pr(\sigma) \cap \pr(\C) = \{(p_1,I_1),\ldots, (p_\ell,I_\ell)\}.$$
In terms of the bipartite graph $G(\C)$ this means that $N(\sigma)= \{(p_1,I_1),\ldots, (p_\ell,I_\ell)\}$. Lastly, we set $\sigma(\C,p_1,\ldots,p_\ell) = \sigma$.
\end{proof}

Applying the construction in the proof of this lemma to Example~\ref{ex:staircases} yields
$$\pi(M_1,5) = \{(6,[1,2]),(5,[3,5])\}\myand \pi(M_2,2) = \{(3,[6,6]), (2,[7,8]), (0,[9,\infty])\},$$
whose union is the profile of the partition $(6,6,5,5,5,3,2,2)$.

\begin{lemma}\label{lem:segments are minimal}
	Let $\C\in \sP$ be such that $\pr(\C)$ is a staircase whose maximal overlapping segments are $M_1,\ldots, M_\ell$.  Then, the $\cl(\C)$-minimal sets in $G(\C)$ are precisely the $\ell$ sets $M_i\subseteq \pr(\C)$.  
\end{lemma}

\begin{proof}
Fix some maximal overlapping segment $M$ of $\pr(\C)$ and write its elements as 
$$(p_1, I_1), (p_2, I_2),\ldots, (p_s, I_s),$$
where $p_i>p_{i+1}$.  Let $a$ and $b$ be such that 
$$\bigcup_{i=1}^s I_i = [a,b].$$  
Further, set $\I=\makeset{i\in [1,s]}{ 1<|I_i|}$ and let $E$ be the set of endpoints   among all the intervals in $M$ (including the symbol $\infty$ if required).  Observe that $E\subseteq [a,b]$ and $|E|= |\I|+1$.

We first demonstrate that $N(M) = \cl(\C)$.  For a contradiction, assume that there exists some $\mu \in \cl(\C)$ so that $\mu\notin N(M)$. Denoting each $p_i$-interval for $\mu$ by $J_i$, this means that for each $(p_i,I_i)\in M$, we must have $J_i\subsetneq I_i$. This immediately says that the interval $J_i$ is such that $|J_i\cap E| \leq 1$ and that $J_i=\emptyset$ whenever $i\notin\I$.  Consequently, 
\begin{equation}\label{eq:too short}
|E\cap  (J_1\cup \cdots \cup J_s)|\leq |\I|.	
\end{equation}
On the other hand, we see that
$$p_1\geq \mu_t\geq p_s\textrm{ if and only if } t\in [a,b].$$
The forward direction follows from the fact that $\pr(\C)$ is a staircase and the reverse direction follows from our definition of a maximal overlapping segment and the meaning of profiles.  So we have
$$E\subseteq [a,b]=J_1\cup \cdots \cup J_s.$$  
But this contradicts (\ref{eq:too short}) above since $|E| = |\I|+1$.  We conclude that $\mu\in N(M)$ as needed.

It remains to show that the sets $M_1,\ldots, M_\ell$ are $\cl(\C)$-minimal and that they are the only $\cl(\C)$-minimal sets.   We can establish both simultaneously by showing that for any $P\subseteq \pr(\C)$ with $N(P) = \cl(\C)$, we have $M_j\subseteq P$ for some $j$.  For a contradiction, assume we have such a $P$ for which $M_j\not\subseteq P$ for all  $1\leq j\leq \ell$.  This means that we can choose some $(q_j, K_j)\in M_j\setminus P$ for each $j$. But then the neighbors of the  partition $\sigma(\C,q_1,\ldots, q_\ell)\in \cl(\C)$ are precisely $(q_1, K_1),\ldots, (q_{\ell}, K_\ell)$.  This contradicts the fact that $N(P) = \cl(\C)$ and completes our proof.
\end{proof}

We are at last in a position to prove our final lemma of this section.  Doing so provides the  justification for (\ref{eq:outline:only staircases survive}) in Section~\ref{sec:GF for P(mu,k)}.

\begin{lemma}\label{lem:only staircases survive}
Let $\C$ be a profile class.  Then
	$$\sum_{P\in \C}(-1)^{|P|} = 
	\begin{cases}
		(-1)^{|\pr(\C)|+\seg(\C)+1} &\text{ if } \pr(\C) \text{ is a staircase}\\
		0 &\text{ otherwise}. 
	\end{cases}$$
 
\end{lemma}

\begin{proof}
We first prove that if $\pr(\C)$ is not a staircase then our alternating sum is 0.  If $\pr(\C)$ is not a staircase, then there exist 
$$(p,[a,b]) \in \pr(\C)\myand (q,[c,d])\in \pr(\C)$$ 
such that either $p> q$ and $c<b$ or $p=q$ and $[a,b]\neq [c,d]$. In the second case, neither $[a,b]$ nor $[c,d]$ is a subset of the other, by definition of a profile, so we can assume (by switching the roles of $[a,b]$ and $[c,d]$ if necessary) that $c< b$. So in any case we have two elements of $\pr(\C)$ such that $p\ge q$ and $c< b$. Lemma~\ref{lem:vanishing profiles} now guarantees that our alternating sum is zero. 

Now consider the case where $\pr(\C)$ is a staircase with $m=\seg(\C)$ maximal overlapping segments: $M_1,\ldots, M_m$.  By Lemma~\ref{lem:segments are minimal} these $m$ sets are precisely the $\cl(\C)$-minimal sets and they partition $\pr(\C)$.  This yields the  calculation 
\begin{align*}
\sum_{P\in \C} (-1)^{|P|} = \sum_{\substack{P\subseteq \cl(\C)\\ N(P) = \pr(\C)}} (-1)^{|P|}=\smashoperator{\sum_{\substack{X\subseteq \pr(\C)\\ N(X) = \cl(\C)}}} (-1)^{|X|}
&=\sum_{\substack{T\subseteq [m]\\T\neq \emptyset}} (-1)^{|T|+1} \left(\sum_{ X\subseteq \pr(\C)\setminus M_T} (-1)^{|X|+|M_T|}\right)\\
&=(-1)^{\seg(\C)+|\pr(\C)|+1},
\end{align*}
where the first and second equalities come from Theorem~\ref{thm:Closures,  Profiles,  and Maxima} and Lemma~\ref{lem:bipartite graph1} respectively;  the third follows from Lemma~\ref{lem:segments are minimal} and an application of Inclusion--Exclusion, where 
$$M_T := \bigcup_{i\in T} M_i;$$
and the last equality follows from (\ref{eq:power set}) and the fact that the $M_i$'s partition $\pr(\C)$.

\end{proof}

\section{Staircases \& Augmented Structures}\label{sec:staircases and augmented structures}
The purpose of this section is to fill in the details of the connection between staircases and augmented structures. In particular we give a proof of (\ref{eq:staircases to aug}) in Section~\ref{sec:GF for P(mu,k)}.

We begin by giving an alternative description of the function $\join_S$ where $S\in \stair(h,k)$.  This alternative description is directly motivated by the following example.
      \begin{example}  Consider the partitions
       $$\alpha = (2,2,2) \myand \beta=(2,2,1,1)$$
       and set $P=\{\alpha,\beta\}$.
If we take $\mu = (4,3,3,1)$, then we see that $$\join_P(\mu)=(6,5,5,1) \join (6,5,4,2)=(6,5,5,4,2,1).$$
      On the other hand, the profile for $P$ is the staircase
      $$S=\{(2, [1,3]), (1, [3,4]), (0, [4,\infty])\}.$$  
Now observe that $\join_P(\mu)$ can be calculated directly in terms of $S$ and $\mu$ by taking
      $$(\underbrace{\mu_1,\ \mu_2,\ \mu_3,\ }_{+2}\underbrace{\ \mu_3,\ \mu_4}_{+1}, \underbrace{\ \mu_4}_{+0})= (6,5,5,4,2,1),$$
where the underbraces indicate how much to add to each term.  
      \end{example}
This example suggests our next lemma.

      \begin{lemma}\label{lem: maxima from profiles}
      Let $P\subseteq \P$ have staircase profile 
      $$S= \{(p_1, [a_1, b_1]),\ldots, (p_{s+1}, [a_{s+1}, \infty])\}$$
      where $k=p_1>p_2>\cdots>p_{s+1}=0$. For any $\mu\in\P$ the partition $\join_S(\mu)$ is the partition  
      $$(\underbrace{\mu_{a_1},\ldots, \mu_{b_1}}_{+p_1}, \underbrace{\mu_{a_2},\ldots, \mu_{b_2}}_{+p_2},\ldots,\underbrace{\mu_{a_{s+1}}, \ldots }_{+p_{s+1}}).$$  The weight of this partition is 
      
      $$\displaystyle \sum_{i=1}^{s+1}\sum_{j=a_i}^{b_i}(\mu_j+ p_i).$$

      \end{lemma}
      
      \begin{proof}  The second assertion follows immediately from the first.  To prove the first assertion, we claim that for any $k\geq p_c> p_d\geq 0$ and any $\sigma, \tau \in P$ (equal or not), where $\sigma$ has a part of size $p_c$ and $\tau$ has a part of size $p_d$,
      $$\min \makeset{(\mu+\sigma)_i}{\sigma_i=p_c} > \max \makeset{(\mu+\tau)_{\ell}}{\tau_{\ell}=p_d}.$$
	  This is because if $\sigma_i=p_c$ and $\tau_{\ell}=p_d$ then by definition of the profile, $i\in [a_c, b_c]$ and $\ell\in [a_d, b_d]$.   By definition of a staircase this implies that $b_c\leq b_d$ and hence $\mu_i\geq \mu_{\ell}$.  Therefore  $\mu_i+p_c> \mu_\ell+p_d$, which verifies our claim.
	      
From this claim we see that if $V_c$ is the set
      $$V_c=\makeset{(\mu+\pi)_i}{\pi\in P, \pi_i=p_c},$$
      then $\min V_c> \max V_d$ when $p_c> p_d$.  Since for each $\pi\in P$ we have $\makeset{i}{\pi_i=p_c}\subseteq [a_c, b_c]$,  the multiplicity of any part in $\join_P(\mu)$ that is also in $V_c$ is its multiplicity in any $\mu+\pi$ such that $\makeset{i}{\pi_i=p_c}=[a_c, b_c]$.   This proves the lemma.
     
      \end{proof}

To motivate our development of the connection between staircases and augmented structures, we start with an example.  

\begin{example}
	Consider the following staircase
$$S = \{(6, [1,2]), (4,[2,4]), (3,[4,4]),(2,[5,6]), (0,[6,\infty])\}$$
and let $\mu=(4^3,3^2,2,1)$ .  Then Lemma~\ref{lem: maxima from profiles} tells us that $|\vee_S(\mu)|$ is
$$\begin{array}{lllll}
(6+\mu_1) &+ (6+\mu_2) &&& \\	
&+ (4+\mu_2) + (4+\mu_3) &+ (4+\mu_4) &&\\
&& +(3+\mu_4) && \\
&&& +(2+\mu_5) &+ (2+\mu_6) \\
&&&&+(0+\mu_6) + (0+\mu_7).
\end{array}$$
Observe that the left-overlapping elements in $S$ correspond to the parts of $\mu$ that are counted more than once.  For the moment, let us ignore such terms.  Doing so yields the sum
$$\begin{array}{llll}
(6+\mu_1) + (6+\mu_2) &&& \\	
& + (4+\mu_3) + (4+\mu_4) &&\\
&& +(2+\mu_5) + (2+\mu_6) &\\
&&& + (0+\mu_7).
\end{array}$$
We encode this sum as the sum of the weights $|\mu|+|\sigma|$, where  $\sigma= (6^2,4^2,2^2)$. Note that the multiplicity of each nonzero part of $\sigma$ is the number of times it appears in the non-overlapping sum.

In order to deal with the overlapping terms in the first sum we introduce the concept of $L$'s.  In particular, consider the $L$'s shaded as shown below where $\mu$ is drawn right-justified and to the immediate left of $\sigma$:
\ytableausetup{boxsize=1em}	
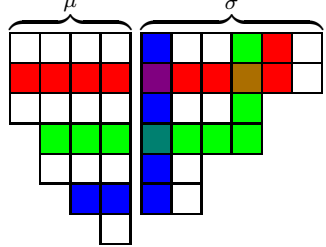
$$
\overbrace{
\ydiagram{4,4,4,1+3,1+3,2+2,3+1}*
[*(red)]{0,4}*
[*(green)]{0,0,0,1+3}*
[*(blue)]{0,0,0,0,0,2+2}
}^{\mu}\, 
\overbrace{
\ydiagram{6,6,4,4,2,2}*
[*(redblue)]{0,1}*
[*(redgreen)]{0,3+1}*
[*(red)]{4+1,1+2}*
[*(red)]{0,4+1}*
[*(greenblue)]{0,0,0,1}*
[*(green)]{3+1,3+1,3+1,1+3}*
[*(blue)]{1,1,1,1,1,1}
}^{\sigma} $$
The columns containing the vertical legs of our $L$'s are chosen as follows.  First take the parts 4, 3, and 0 as they correspond to the left-overlapping intervals in $S$.  Then add one to each to obtain the set $A=\{5,4,1\}$ which will be the set of indices of the columns  containing the vertical legs of our $L$'s.  

In this case, the vertical parts of the blue, green, and red $L$'s contain $\sigma^*_1$, $\sigma^*_4$, and $\sigma^*_5$ boxes, respectively.    The key observation is that we can now write our first sum in terms of these $L$'s.  Consider
\begin{align*}
|\vee_S(\mu)| &= |\mu| + |\sigma| + (4+\mu_2) + (3+\mu_4) + (0+\mu_6)\\
&= |\mu| + \underbrace{\sigma_2^*+\sigma_3^*+\sigma_6^*}_{\text{columns $\notin A$}} +  \underbrace{\sigma^*_1 + \sigma^*_4 + \sigma^*_5}_{\text{vertical parts of $L$'s}} + \underbrace{ (4+\mu_{\sigma^*_5})+ (3+ \mu_{\sigma^*_4})+( 0 + \mu_{\sigma^*_1} )}_{\text{horizontal parts of $L$'s (excluding corners)}}.
\end{align*}
\end{example}

In light of this example, we make the following definition.
\begin{definition}
We say an element $(\sigma,A)\in \P(h,k)\times 2^{[1,k]}$ is a \emph{marked partition} since we think of the set $A$ as indicating the marked columns of $\sigma$. We let 
$$\M(h,k)=\makeset{(\sigma,A)\in \P(h,k)\times 2^{[1,k]}}{\sigma^*_i>1\text{ for all } i\in A}.$$  	
\end{definition}

\begin{lemma}\label{lem:bijection stair to marked columns}
	Fix $h,k > 0$.  There exists a bijection $f:\stair(h,k) \to \M(h,k)$ so that if $f(S)=(\sigma,A)$ and $\mu$ is any partition, then
	\begin{enumerate}[label=(\roman*)]
		\item $|A| = |S|-\seg(S)$, and
		\item $\displaystyle |\join_S (\mu)| = |\mu| + |\sigma|+ \sum_{i\in A} ( \mu_{\sigma^*_i}+i-1)$.
	\end{enumerate}  
\end{lemma}
\begin{proof}
	Fix $S\in \stair(h,k)$ and denote its elements as 
	$$S = \{(p_1,[a_1,b_1]), (p_2,[a_2,b_2]),\ldots, (p_s,[a_s,b_s]),(p_{s+1},[a_{s+1},b_{s+1}])\},$$
	where $k=p_1>p_2>\cdots > p_s>p_{s+1}=0$ and $a_2\neq 1$. Define $I$ to be the set of indices that correspond to left-overlapping elements in $S$.  Let
	$$A = \makeset{1+p_i}{i\in I}\myand \sigma = (p_1^{m_1},p_2^{m_2},\cdots, p_{s}^{m_{s}}),$$ 
	where 
$$	m_i = \begin{cases}
		b_i - a_i +1&\text{ if } i\notin I\\
		b_i - a_i&\text{ if } i\in I.
	\end{cases}$$
Set $f(S) = (\sigma, A)$.  With these definitions we see that 

\begin{equation}\label{eq:a_i}
a_i = \begin{cases}
m_1+\cdots+m_{i-1}+1 & \text{ if } i\notin I\\
m_1+\cdots+m_{i-1} &\text{ if } i\in I
\end{cases}	
\end{equation}
for $1\leq i\leq s+1$. We also observe that any $i\in I$ must be at least 2.  Consequently Property \emph{(ii)} in our definition of staircases implies that, for $i\in I$,
\begin{equation}\label{eq:overlapping endpoint}
2\leq a_i = m_1+\cdots +m_{i-1} = \sigma^*_{1+p_i}.
\end{equation}

With our definitions and basic observations established our immediate task is to show that $f$ is well defined.  To this end, we first show $\sigma\in \P(h,k)$.  To verify this, observe that  $m_1+\cdots + m_s$ is the length of $S$, so $h_\sigma< h$.  Additionally, $w_\sigma = k$ as $p_1=k$ and $m_1\neq  0$. So $\sigma\in  \P(h,k)$.  It is immediate that $A\subseteq [1,k]$ since $1\notin I$.  Lastly,  (\ref{eq:overlapping endpoint}) guarantees that $(\sigma, A) \in \M(h,k)$.  

Our next task is to show $f$ is bijective.  We first note that we can easily recover the parts $p_1,\ldots, p_s$ from the pair $(\sigma, A)$.  In turn, knowing these parts and the set $A$ allows us to recover the index set $I$.  Using $\sigma$, $I$, and (\ref{eq:a_i}) we can recover the left endpoints $a_i$ of our staircase.  Finally, the left endpoints together with $I$ allow us to reconstruct the right endpoints $b_i$.  This inverse construction proves that $f$ is 1-1.  Furthermore, let us consider applying this inverse construction to any $(\tau,B)\in\M(h,k)$ to obtain 
$$T = \{(q_1,[c_1,d_1]), (q_2,[c_2,d_2]),\ldots, (q_t,[c_t,d_t]), (q_{t+1},[c_{t+1},\infty])\},$$
where $k=q_1>q_2>\cdots > q_t>q_{t+1}=0$.  Provided that $c_2\geq 2$ it easily follows that $T\in\stair(h,k)$.   To address this proviso, we see from (\ref{eq:a_i}) that if $(q_2,[c_2,d_2])$ is not left-overlapping then $c_2\geq 2$.  On the other hand, the only way this element can be left-overlapping is if $q_2+1\in B$, in which case our definition of $\M(h,k)$ tells us that $1<\tau_{q_2+1}^*$. This says that there are at least two parts in $\tau$ (counting multiplicities) that are greater than $q_2$.  As the only part greater than $q_2$ is $q_1$ we conclude that $1<d_1\leq c_2$ as needed. 

We now shift our attention to Properties \emph{(i)} and \emph{(ii)}. It is clear from our definitions of $I$ and $A$ that Property \emph{(i)} holds. Property \emph{(ii)} is shown to hold via the following computation.   By Lemma~\ref{lem: maxima from profiles}, we have
\begin{align*}
|\join_S(\mu)| &=\sum_{i=1}^{s+1}\sum_{j= a_i}^{b_i}(\mu_j+ p_i)\\
&=\sum_{i\not\in I}\sum_{j=a_i}^{b_i}(\mu_j+ p_i) + \sum_{i\in I}\sum_{j=a_i+1}^{b_i}(\mu_j+ p_i)+ \sum_{i\in I}(\mu_{a_i}+ p_i)\\
&= |\mu| + |\sigma| + \sum_{i\in I}(\mu_{a_i}+ p_i),
\shortintertext{which by (\ref{eq:overlapping endpoint}) and our definition of $A$,}
&= |\mu| + |\sigma| + \sum_{i\in A} (\mu_{\sigma^*_i} + i-1).
\end{align*}
\end{proof}

Lastly, we need to bijectively transform marked partitions into augmented structures.  To do this we define, for any $\mu$, the mapping 
$$g:\M(h,k) \to \aux(\mu,h,k)$$
by setting $g(\sigma, A) = (\mu,\lambda, \off)$ where the columns of $\lambda$ and $\off$ are given by the multisets
$$\makeset{\sigma^*_i}{i\in A}\myand \makeset{|A_{>i}|+\sigma^*_i}{i\not\in A},$$
respectively.  To explain this definition, let us return to the above example of a marked partition: 
$$
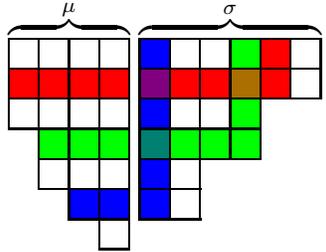
\overbrace{
\ydiagram{4,4,4,1+3,1+3,2+2,3+1}*
[*(red)]{0,4}*
[*(green)]{0,0,0,1+3}*
[*(blue)]{0,0,0,0,0,2+2}
}^{\mu}\, 
\overbrace{
\ydiagram{6,6,4,4,2,2}*
[*(redblue)]{0,1}*
[*(redgreen)]{0,3+1}*
[*(red)]{4+1,1+2}*
[*(red)]{0,4+1}*
[*(greenblue)]{0,0,0,1}*
[*(green)]{3+1,3+1,3+1,1+3}*
[*(blue)]{1,1,1,1,1,1}
}^{\sigma}.$$
If we ``push" all the unmarked columns past and to the right of the marked ones we obtain the configuration
$$\overbrace{
\ydiagram{4,4,4,1+3,1+3,2+2,3+1}*
[*(red)]{0,4}*
[*(green)]{0,0,0,1+3}*
[*(blue)]{0,0,0,0,0,2+2}
}^{\mu}\,\overbrace{\ydiagram{3,3,2,2,1,1}*
[*(redblue)]{0,1}*
[*(redgreen)]{0,1+1}*
[*(red)]{2+1,2+1}*
[*(greenblue)]{0,0,0,1}*
[*(green)]{1+1,1+1,1+1,1+1}*
[*(blue)]{1,1,1,1,1,1}}^{\lambda}\quad  \overbrace{\ydiagram{3,3,2,2,1,1}}^{\off'}.$$ 
Observe that this transformation shortened the horizontal legs of the green and red $L$'s by exactly the number of unmarked columns to their left. This discrepancy is accounted for by defining the columns of $\off$ to be $|A_{>i}|+\sigma^*_i$ for each $i\notin A$.  

  To establish that this mapping is bijective we will prove that for any $C>0$ there is a unique $0\leq i\leq w_\lambda$  so that placing a column of length $C-i$  immediately to the left of the $i$ rightmost columns of $\lambda$ yields a partition.  For example, if $C = 5$ and $\lambda$ is as above, then among the four possible configurations  
  $$
  \begin{tabular}{cccc}
    \ydiagram{3,3,2,2,1,1}*[*(lightgray)]{3+1,3+1,3+1,3+1,3+1}\hspace{5ex}   &
    \ydiagram{2,2,2,2,1,1}*[*(lightgray)]{2+1,2+1,2+1,2+1}*[*(white)]{3+1,3+1}\hspace{5ex} &
    \ydiagram{1,1,1,1,1,1}*[*(lightgray)]{1+1,1+1,1+1}*[*(white)]{2+2,2+2,2+1,2+1}\hspace{5ex} &
    \ydiagram[*(lightgray)]{1,1,1}*[*(white)]{1+3,1+3,1+2,1+2,1+1,1+1}\hspace{5ex}\\
  \end{tabular}
  $$
  only the second, when $i=1$, yields a partition.  The following lemma proves this always occurs.

\begin{lemma}
Fix any partition $\lambda$  and denote its column lengths (from right to left) by $c_1\leq c_2\leq c_3\leq \cdots \leq c_{w_{\lambda}}$.   For any $C>0$ there is a unique $0\leq i\leq w_\lambda$ so that 
\begin{equation}\label{eq:unique insertion}
c_{i+1}\geq C-i \geq c_{i},
\end{equation}
where we set $c_{w_\lambda+1} =\infty$ and $c_0 = 0$.  
\end{lemma}
\begin{proof}
	Define $i = \min\makeset{j\in[0,w_\lambda]}{c_{j+1} \geq C-j }$. (This set is not empty, and hence $i$ is well defined, since $c_{w_\lambda+1}=\infty$.) By definition of $i$ it follows that if $i> 0$ then $C-(i-1) > c_i$.  Hence (\ref{eq:unique insertion}) holds in any case.	Uniqueness follows since we cannot have indices $j>i$ such that 
$$c_{j+1}\geq C-j \geq c_{j}\geq c_{i+1} \geq C-i \geq c_{i},$$
as then $j\leq i$.   
\end{proof}

Lastly, if $S\in \stair(h,k)$ and $f(S)= (\sigma,A)$ and $g(\sigma, A) = (\mu,\lambda, \off)$, then we claim that, for all $\mu$, $$|\join_S(\mu)|= |(\mu,\lambda,\off)|.$$
In fact this follows via the following computation:
\begin{align*}
	|\join_S(\mu)| &=  |\mu| + |\sigma|+ \sum_{i\in A} (\mu_{\sigma^*_i}+i-1),\\ 
	\intertext{which, by segregating the $i-1$ columns to the left of column $i$ into marked and unmarked ones,} 
	&= |\mu| + \sum_{j\in B} \sigma^*_j + \sum_{i\in A} (\sigma^*_i+ \mu_{\sigma^*_i} + |A_{<i}| + |B_{<i}|),
\intertext{where $B = [1,k]\setminus A$. By noting that the marked columns in $\sigma$ are precisely the columns of $\lambda$, this}
	&=|\mu| + \sum_{j\in B} \sigma_j^* +\frac{a(a-1)}{2}+  \sum_{\ell=1}^a \left( \lambda^*_\ell +\mu_{\lambda^*_\ell} \right) + \sum_{i\in A} |B_{<i}|,\\
		\intertext{where $a =|A|$, which by the short argument below, }
	&=|\mu| + \sum_{j\in B} \left(\sigma_j^* + |A_{>j}|\right) +  \frac{a(a-1)}{2}+\sum_{\ell=1}^a \mu_{\lambda^*_\ell} + \lambda^*_\ell\\
	&= |\mu| + |\off| + |\lambda| + \frac{a(a-1)}{2} + \sum_{\ell=1}^a \mu_{\lambda^*_\ell}\\
	&= |(\mu,\lambda,\off)|.
\end{align*}
The following calculation justifies the fourth equality above:
$$
\sum_{i\in A} |B_{<i}| = \sum_{i\in A}\sum_{j\in B} 1_{\{j<i\}} = \sum_{j\in B}\sum_{i\in A} 1_{\{j<i\}}	=  \sum_{j\in B} |A_{>j}|,	
$$
where $1_{\{\cdot\}}$ is the usual indicator function.  

We summarize our findings from this section with a final lemma.
\begin{lemma}\label{lem:bijection between stair and aux}
	For any $\mu$, the composition of maps
	 $$g\circ f: \stair(h,k)\to \aux(\mu,h,k)$$ 
	 is a bijection such that if $S\mapsto (\mu,\lambda,\off)$, then
	\begin{enumerate}[label=(\roman*)]
		\item $w_{\lambda} = |S|-\seg(S)$, and
		\item $|\join_S (\mu)| = |(\mu,\lambda,\off)|$.
	\end{enumerate}  

\end{lemma}

\begin{ack}

The first version of this work, written before we were aware of \cite{foata1970rook}, was formulated in terms of Wilf equivalence and a condition involving a simplified version of the $(i,j)$-transforms of Foata and Sch\"{u}tzenberger.  We thank Bruce Sagan  for pointing out to us the work of Foata and Sch\"{u}tzenberger in \cite{foata1970rook}, and the connection between our original condition and rook equivalence.
\end{ack}

\bibliography{mybib}
\bibliographystyle{siam}      
    \end{document}